\documentclass{amsproc}

 \newtheoremstyle{plain}
     {3pt}
     {3pt}
     {}
     {}
     {\sc}
     {:}
     {.5em}
     {}

 \newtheoremstyle{tshtheorem}
     {3pt}
     {5pt}
     {\it}
     {}
     {\sc}
     {.}
     {.5em}
     {}

\theoremstyle{tshtheorem}
\newtheorem{theorem}{Theorem}[section]

\newtheorem{thmd}[theorem]{Darboux's Theorem}

\newtheorem{thmc}[theorem]{Convexity Theorem}
\newtheorem{thmi}[theorem]{Injectivity Theorem}
\newtheorem{thms}[theorem]{Surjectivity Theorem}

\newtheorem{thml1}[theorem]{First Localization Theorem}
\newtheorem{thml2}[theorem]{Second Localization Theorem}
\newtheorem{thml3}[theorem]{Third Localization Theorem}

\theoremstyle{plain}

\newtheorem{remark}[theorem]{Remark}
\newtheorem{assum}{Assumption}

\numberwithin{equation}{section}
\numberwithin{figure}{section}

\usepackage{epsfig}
\usepackage{psfrag} 
\input xy
\xyoption{all}
\usepackage{amsmath}
\usepackage{amssymb}
\usepackage{amscd}

\newcommand{\C}{{\mathbb{C}}}

\newcommand{\Z}{{\mathbb{Z}}}
\newcommand{\Q}{{\mathbb{Q}}}
\newcommand{\R}{{\mathbb{R}}}

\newcommand{\X}{{\mathcal{X}}}
\newcommand{\Y}{{\mathcal{Y}}}

\newcommand{\F}{{\mathcal{F}}}

\newcommand{\algt}{\mathfrak{t}}
\newcommand{\algg}{\mathfrak{g}}

\newcommand{\Id}{\mathrm{Id}}

\begin{document}

\title[Act globally, compute locally]{Act globally, compute locally:\\ group actions, fixed points, and localization}

\author{Tara S. Holm}
\address{Department of Mathematics, Cornell University, Ithaca, NY 14853-4201 USA}
\email{tsh@math.cornell.edu}
\thanks{TSH is grateful for the support of the NSF through the grant DMS-0604807.}

\date{April 2007}


\begin{abstract}

Localization is a topological technique that allows us to make global equivariant computations in terms of local data at the fixed points. For example, we may compute a global integral by summing integrals at each of the fixed points. Or, if we know that the global integral is zero, we conclude that the sum of the local integrals is zero. This often turns topological questions into combinatorial ones and vice versa.
This expository article features several instances of localization that occur at the crossroads of symplectic and algebraic geometry on the one hand, and combinatorics and representation theory on the other.  The examples come largely from the symplectic category, with particular attention to toric varieties. In the spirit of the 2006 International Conference on Toric Topology at Osaka City University, the main goal of this exposition is to exhibit toric techniques that arise in symplectic geometry.
\end{abstract}

\maketitle

\tableofcontents

\section*{Introduction}
 
In topology, a {\bf localization result} relates the $G$-equivariant topology of $X$ to the topology of the $G$-fixed points and their local isotropy data.  Results of this nature go back to the fundamental work of Borel; this will be the starting point for our discussion of localization.  We then turn to work of Atiyah-Bott and Berline-Vergne; and finally of Guillemin-Ginzburg-Karshon.  Our motivation is to place these results in the context of symplectic geometry, and to give an overview of the ways in which localization has played a role in recent research in symplectic geometry.

To this end, the article is organized as follows.  Section~\ref{sec:symp} provides a quick introduction to symplectic geometry, including a description of the symplectic point of view on toric varieties.   We turn to several localization results in Section~\ref{sec:local}, and describe how Borel's result simplifies topological computations (in the symplectic category) in Section~\ref{sec:compute}.  We conclude with a combinatorial version of localization in Section~\ref{sec:polytopes}.  We have included a (by no means exhaustive) list of references where the reader may find additional details.

\medskip
\noindent {\bf Acknowledgments:}  This article grew out of two lectures, one at the 2005 Summer Institute on Algebraic Geometry at the University of Washington, and the other at the 2006 International Conference on Toric Topology at Osaka City University.  Both conferences were extremely fruitful for the communication of mathematics, and I am very grateful for the opportunity to participate in each.  The organizers of the Osaka conference, Megumi Harada, Yael Karshon, Mikiya Masuda, and Taras Panov, deserve particular mention for making a daunting visit to Japan a very well-organized and enjoyable one.  I extend many thanks to the referee and to Megumi Harada for very carefully reading an earlier draft of this paper, and providing extensive comments.

\section{A brief review of the symplectic category}\label{sec:symp}

We begin with a very quick introduction to symplectic geometry.  Many more details can be found in \cite{CdS:book, mcduff-salamon}.
Let $M$ be a {\bf manifold} with a {\bf symplectic form} on $M$, that is, a two-form 
$\omega\in \Omega^2(M)$ that is closed ($d\omega = 0$) and non-degenerate (the top power of $\omega$ is a volume form on $M$).  
In particular, the non-degeneracy condition implies that $M$ must be an even-dimensional manifold.  The key examples include
\begin{enumerate}
\item $M=S^2 = \C P^1$ with $\omega_p (\X,\Y)=$ signed area of the parallelogram spanned 
by $\X$ and $\Y$;

\item $M$ any Riemann surface with $\omega$ as in (1);

\item $M= \R^{2d}$ with $\omega = \sum dx_i\wedge dy_i$; and

\item $M= \mathcal{O}_\lambda$ a coadjoint orbit of a compact connected semisimple Lie group $G$, 
equipped with the Kostant-Kirillov-Souriau form $\omega$.  For the group $G=SU(n+1)$, this class of examples includes complex projective space $\C P^n$, the full flag variety $\F \ell (\C^{n+1})$, and all other partial flag varieties.
\end{enumerate}
Example (3) gains particular importance because of 
\begin{thmd}
Let $M$ be a $2d$-dimensional symplectic manifold with symplectic form $\omega$.  
Then for every point $p\in M$, there exists a coordinate chart $U$ about 
$p$ with coordinates $x_1,\dots,x_d,y_1,\dots,y_d$ so that on this chart,
$$
\omega = \sum_{i=1}^d dx_i\wedge dy_i.
$$
\end{thmd}
\noindent Thus, whereas Riemannian geometry uses local invariants such as curvature 
to distinguish metrics, symplectic forms are locally indistinguishable.

The symmetries of a symplectic manifold may be encoded by a group action.  
Here we restrict ourselves to a compact connected abelian group $T= (S^1)^n$.  
An action of $T$ on $M$ is {\bf symplectic} if it preserves $\omega$; that is, 
$\rho_g^*\omega = \omega$, for each $g\in T$, where $\rho_g$ is the diffeomorphism 
corresponding to the group element $g$.  The action is {\bf Hamiltonian} if in 
addition, for every vector $\xi\in\algt$ in the Lie algebra $\algt$ of $T$, the vector field
$$
\X_\xi(p) = \frac{d}{dt}\Big[ \exp (t\xi)\cdot p \Big] \bigg|_{t=0}
$$
is a {\bf Hamiltonian vector field}.  That is, we require that 
\begin{equation}\label{eq:mmap}
\omega(\X_\xi, \cdot ) = d\phi^\xi
\end{equation}
is an exact one-form.  Thus each $\phi^\xi$ is a smooth function on $M$ defined by the differential equation \eqref{eq:mmap}, so determined 
up to a constant.  Taking them together, we may define a {\bf moment map}\footnote{The map $\Phi$ is also called a {\bf momentum} map.}
$$
\begin{array}{rcc}
\Phi: M  & \to & \algt^* \\
 p & \mapsto & \left(\begin{array}{rcl}
                \algt & \longrightarrow & \R \\
                \xi & \mapsto & \phi^\xi(p)
                \end{array}\right).
\end{array}
$$

We now examine how substantial the Hamiltonian assumption is.  Certainly, if $M$ is compact, then any component of the moment map $\phi^\xi$  has bounded image.  Thus, as a smooth function, $\phi^\xi$ must have critical points, points where $d\phi^\xi = 0$.  By the moment map condition \eqref{eq:mmap}, this means $\omega(\X_\xi, \cdot ) = 0$.  The form $\omega$ is non-degenerate, so the vector field $X_\xi$ must be $0$ at the critical points, which implies that  the subgroup $\exp(t\xi)$ acts trivially on these critical points: they are fixed.  Thus, a Hamiltonian circle action on a compact symplectic manifold necessarily has fixed points. 

Frankel shows that in the {\em K\"ahler} setting, the existence of fixed points is a sufficient condition for a circle action to be Hamiltonian.  For every symplectic form on a manifold, there is a compatible {\bf almost complex structure} $J: TM\to TM$ with $J^2 = -\Id$  When this almost complex structure is {\bf integrable}, then we say that the manifold is {\bf K\"ahler}.

\begin{theorem}[Frankel \cite{frankel}]
A symplectic circle action which preserves the  compatible complex structure on a compact K\"ahler manifold $M$ is Hamiltonian if and only if it has fixed points.
\end{theorem}

The K\"ahler hypothesis is strictly necessary; indeed McDuff has addressed the question for general symplectic manifolds.

\begin{theorem}[McDuff \cite{mcduff}]\label{thm:mcduff}
$\phantom{.}$
\begin{enumerate}
\item[(a)] If $M$ is a compact four-dimensional symplectic manifold, a symplectic circle action on $M$ is Hamiltonian if and only if it has fixed points.
\item[(b)] There is a symplectic circle action on a symplectic six-dimensional manifold $M$ that has fixed points, but is not Hamiltonian.
\end{enumerate}
\end{theorem}

Returning to our examples, we have Hamiltonian actions in all but the second example.

\begin{enumerate}

\item The circle $S^1$ acts on $M=S^2 = \C P^1$ by rotations, with the north and south
poles fixed.  If we use angle $\theta$ and height $h$ coordinates on $S^2$, then the vector 
fields this action generates are tangent to the latitude lines.  Identifying $\algt\cong \R$ and choosing $\xi = 1\in\R$,
the vector field is
$\X_\xi = \frac{\partial}{\partial\theta}$. Since the symplectic form $\omega = d\theta\wedge dh$, 
the moment map condition \eqref{eq:mmap} is $\omega(\X_\xi, \cdot ) = dh$, and a moment map is the height function, 
shown in Figure~\ref{fig:sphere}.
\begin{figure}[h]
  \begin{center}
    \epsfig{file=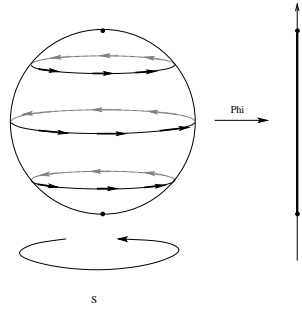,height=4cm}
    \caption{The vector field and moment map for $S^1$ acting by rotations on $S^2$.}
    \label{fig:sphere}
  \end{center}
 \end{figure}

\item If $M$ is a two-torus $M=T^2 = S^1\times S^1$, then $S^1\times S^1$ acts 
on itself by multiplication.  This action is certainly symplectic.  It is a free action, though, and as Hamiltonian actions necessarily have fixed points, it is not Hamiltonian.    
In fact, no Riemann surface of non-zero genus admits a nontrivial Hamiltonian torus action.

\item The torus $T^d$ acts on $M= \R^{2d} = \C^d$ by rotation of each copy of $\C = \R^2$.  
This action is Hamiltonian, and identifying $\algt^* \cong \R^d$, a moment map is
$$
\Phi (z_1,\dots,z_n) = (|z_1|^2, \dots,|z_d|^2).
$$
The moment map does depend on the coordinates that we choose on the torus, so depending on whether an author identifies $T^d \cong \R^d/\Z^d$ or $T^d\cong \R^d/ 2\pi\Z^d$, a moment map may have an extra factor of $2\pi$.

\item Each coadjoint orbit $M= \mathcal{O}_\lambda\subseteq \algg^*$ may be identified 
as a homogeneous space $G/L$, where $L$ is a Levi subgroup of the Lie group $G$.  Thus $G$, 
and hence the maximal torus $T$, acts on $M$ by left multiplication.  A $G$-moment map is 
inclusion
$$
\Phi_G : \mathcal{O}_\lambda\hookrightarrow \algg^*,
$$
and a $T$-moment map is the $G$-moment map composed with the natural projection 
$\algg^* \to\algt^*$ that is dual to the inclusion $\algt\hookrightarrow\algg$. 
\end{enumerate}

A {\bf localization} phenomenon is a global feature of the $T$ action on $M$ that  can be completely determined or described by the evidence of that feature at the $T$-fixed points.  An example of this type of occurrence is the following theorem, due independently to Atiyah and to Guillemin and Sternberg.

\begin{thmc}[Atiyah \cite{At:convexity}, Guillemin-Sternberg \cite{GS:convexity}]\label{thm:convex}
If $M$ is a compact Hamiltonian $T$-space, then $\Phi(M)$ is a convex polytope.  It is the convex hull of  the images $\Phi(M^T)$ of the $T$-fixed points.
\end{thmc}

In proving this theorem, Atiyah also establishes the fact that any component $\phi^\xi$ of the moment map is a {\bf Morse function} on $M$ (in the sense of Bott), with critical set the fixed set of the subgroup of elements $\exp (t\xi)$ for $t\in \R$.  For almost all $\xi$, this fixed set is precisely the set of torus fixed points $M^T$.  This is a strong indication that the topology of $M$ is dictated by the topology of the fixed point set.  We will return to this theme in the next section.

For an {\bf effective}\footnote{An action is effective if no positive dimensional subgroup acts trivially.} Hamiltonian $T$ action on $M$, 
$$
\dim(T)\leq \frac{1}{2}\dim(M).
$$
We say that the action is {\bf toric} if this inequality is in fact an equality.  A symplectic manifold $M$ with a toric Hamiltonian $T$ action is called a {\bf symplectic toric manifold}.  There is not consensus in the literature about this term: some authors use {\bf toric variety} to parallel the topic in algebraic geometry.  Every symplectic toric manifold is a toric variety, in the sense of  \cite{fulton}, but the converse is not true (see \cite[p.\ 25--26]{fulton}).  In this paper, we remain firmly grounded in the symplectic category and stick to the terminology symplectic toric manifold.  

Delzant used the moment polytope to classify symplectic toric manifolds.  
To understand this result, we must first review several features of convex polytopes; a more comprehensive introduction to polytopes can be found in \cite{Zieg}.  A polytope $\Delta$ in $\R^n$ may be defined as the convex hull of a set of points, or alternatively as a (bounded) intersection of a finite number of half-spaces in $\R^n$.  
We say $\Delta$ is {\bf simple} if there are $n$ {\bf edges} adjacent to each {\bf vertex}, and it is {\bf rational} if the edges have rational slope. For a vector with rational slope, the {\bf primitive vector} with that slope is the shortest positive multiple of the vector that is in the lattice $\Z^n\subseteq \R^n$. A  simple polytope is {\bf smooth at a vertex} if the $n$ primitive edge vectors emanating from the vertex span the lattice $\Z^n\subseteq\R^n$ over $\Z$. It is {\bf smooth} if it is smooth at each vertex.  We may now state Delzant's result.

\begin{theorem}[Delzant \cite{Del}]
There is a one-to-one correspondence
$$
\left\{\begin{array}{c}
\mbox{compact symplectic}\\
\mbox{toric manifolds}
\end{array}\right\}
\leftrightsquigarrow
\left\{\begin{array}{c}
\mbox{simple rational}\\
\mbox{smooth convex polytopes}
\end{array}\right\} .
$$
\end{theorem}

To each symplectic toric manifold, the polytope we associate to it is its moment polytope.  It is not hard to check that the moment polytope has the aforementioned properties.  To such a
polytope, on the other hand, we may construct a symplectic toric manifold via {\em symplectic reduction}.
The adjective {\bf compact} in Delzant's theorem can be relaxed.  To do so, we must replace {\bf polytopes} with {\bf polyhedra}, which are possibly unbounded.  The other conditions on the polyhedra remain.  For further details, the reader should consult \cite[Chapter VII]{audin}.

The moment map $\Phi:M\to\algt^*$ is a $T$-invariant map: it maps entire $T$-orbits to the same point in $\algt^*$.  
Thus, if $\alpha\in\algt^*$ is a regular value, then the level set $\Phi^{-1}(\alpha)$ is a 
$T$-invariant submanifold of $M$.    Moreover, the action of $T$ on the level set is {\bf locally free}: 
it has only finite stabilizers. We deduce this fact  by examining the moment map condition. If $\alpha$ is a critical value for $\Phi$, then it is critical for 
some $\phi^\xi$.  The moment map condition \eqref{eq:mmap}, together with the 
non-degeneracy of $\omega$, then implies that the vector field $\X_\xi$ must be zero, which means that the critical points must be fixed by the subgroup $\exp(t\xi)$.  Thus, when $\alpha$ is a regular value, it is not critical for any $\phi^\xi$, and so {\bf no} positive-dimensional subgroup of $T$ fixes any point in the level set.  We deduce, then, that $T$ acts locally freely on $\Phi^{-1}(\alpha)$.  Therefore, the 
{\bf symplectic reduction} $M/\!/T(\alpha) = \Phi^{-1}(\alpha)/T$ is an orbifold.  In fact,

\begin{theorem}[Marsden-Weinstein \cite{MW:reduction}]
If $M$ is a Hamiltonian $T$-space and $\alpha$ is a regular value of the moment map $\Phi$, then 
the symplectic reduction $M/\!/T(\alpha)$ is a {\bf symplectic} orbifold. 
\end{theorem}

Symplectic reduction is an important method of constructing new symplectic manifolds 
from old. From our examples, we may construct several classes of symplectic manifolds, including symplectic toric manifolds.

\begin{enumerate}

\item $S^1$ acts on $M=S^2 = \C P^1$ by rotations.  The level set of a regular value is a 
latitude line, which the circle rotates.  The quotient is a point.

\item $T^d$ acts on $M= \R^{2d} = \C^n$ by rotation of each copy of $\C = \R^2$.  
The level set of a regular value is a copy of $T^d$, and so again the quotient is a point. 
We may also restrict our attention to a subtorus $K\subseteq T^n$.  
The action of $K$ is still Hamiltonian, and for 
certain choices of $K$ and $\alpha$, $\C^n/\!/K (\alpha)$ is a {\bf symplectic toric manifold} with the residual $T^n/K$ action.  Given a rational simple smooth convex polytope $\Delta$, Delzant found
a judicious way to choose $K$, using the combinatorics of $\Delta$, so that
$\C^n/\!/K (\alpha)$ has moment polytope precisely $\Delta$.

\begin{figure}[h]
  \begin{center}
  \psfrag{CC}{$\C^3$}
  \psfrag{Phi}{$\Phi_T$}
  \psfrag{Pi}{$\Phi_K$}
 \psfrag{aa}{$\alpha$}
      \epsfig{file=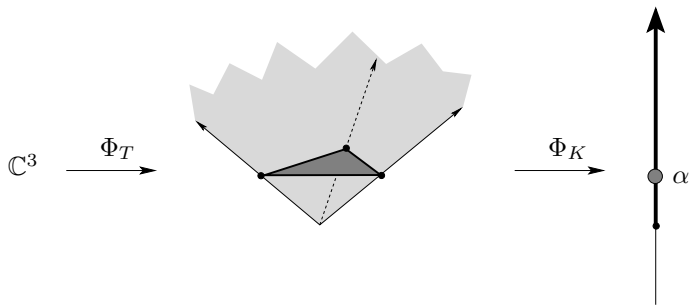,height=4cm}
    \caption{The moment map for $T=T^3$ and $K=S^1$ acting on $\C^3$.  The triangle represents the moment polytope for the reduction $\C^3/\! /S^1(\alpha)$.}
  \end{center}
 \end{figure}

\end{enumerate}

\begin{figure}[h]
  \begin{center}
  \psfrag{=}{$=$}
  \psfrag{-}{$-$}
  \psfrag{+}{$+$}
      \epsfig{file=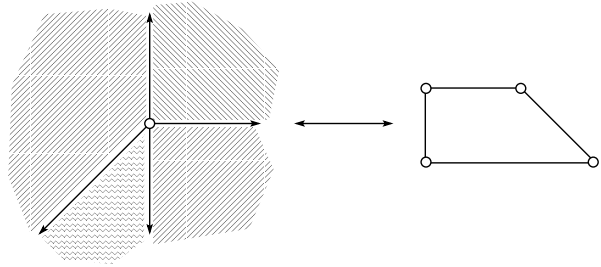,height=3.5cm}
    \caption{The fan and  corresponding polytope for a Hirzebruch surface.}
    \label{fig:fan}
  \end{center}
 \end{figure}

We conclude this section with a brief comment on the relationship between symplectic toric manifolds  and toric varieties  in the algebraic geometry sense.  In the latter case, the combinatorial data used to construct a toric variety is a {\bf fan}.  Given a polytope, the corresponding fan is generated by the inward pointing normal vectors to the facets.  There is a family of polytopes (dilations of one another) corresponding to any fan.  Choosing a particular polytope corresponds to specifying an invariant ample line bundle on the (algebraic) toric variety.  An example is shown in Figure~\ref{fig:fan}.  
Toric varieties in the algebraic category may be constructed using {\bf geometric invariant theory} quotients.  It is possible to prove that these two very different prescriptions yield homeomorphic quotients.  There are many further details on fans in \cite{fulton} and on geometric invariant theory in \cite{mfk}.

\section{Equivariant cohomology and localization theorems}\label{sec:local}

We have seen that a generic component $\phi^\xi$ of the moment map is a Morse function, in the sense of Bott, on a Hamiltonian $T$-manifold $M$, with critical set $M^T$.  This is a certainly a form of localization, and it allows us to compute bounds for the Betti numbers of $M$ from those of $M^T$. It does not, however, give us any information about the ring structure of $H^*(M;\Q)$.  We can use a stronger form of localization to compute the {\bf equivariant cohomology} of the $T$-manifold, including its ring structure.  In fact, from the equivariant cohomology ring, we may deduce results about ordinary cohomology ring as well.

Equivariant cohomology is a generalized cohomology theory in the equivariant category.  We use the Borel model to compute equivariant cohomology.  For a topological group $G$, let $EG$ be a contractible space on which $G$ acts freely (such a space always exists, though it may be infinite-dimensional).  Then for any $G$-space $X$, the diagonal action of $G$ on $X\times EG$ is free, and 
 $$
 X_G = (X\times EG)/G
 $$
is the {\bf Borel mixing space} or {\bf homotopy quotient} of $X$.  We define the equivariant cohomology ring to be
 $$
 H_G^*(X;R) := H^*(X_G;R),
 $$
 where $H^*(-;R)$ denotes singular cohomology with coefficients in the commutative ring $R$.  Thus, when $X$ is a free $G$-space, we may identify
 $$
 H_G^*(X;R) \cong H^*(X/G;R).
$$
At the other extreme, if $G$ acts trivially on $X$, then
 $$
 H_G^*(X;R) \cong H^*(X\times BG;R),
$$
where $BG = EG/G$ is the {\bf classifying space} of $G$.  Note that the cohomology of the classifying space, $H^*(BG;R) \cong H_G^*(pt;R)$, is the equivariant cohomology ring of a point.

\begin{remark}
Cartan gave an alternative description of equivariant cohomology in terms of equivariant differential forms, in the case when $X$ is a manifold.  Full details, including a proof of that Cartan's description is equivalent to the Borel model, may be found in \cite{GS:super}.
\end{remark}

For a compact torus $T=(S^1)^n$, we may take $EG = (S^\infty)^n$, where $S^\infty\subseteq \C^\infty$ is the (contractible!) infinite dimensional sphere.  Then the classifying space $$BT = ET/T = (\C P^\infty)^n$$ is $n$ copies of infinite complex projective space, whence $H_T^*(pt;\Z)=\Z[x_1,\dots,x_n]$ with $\deg(x_i)=2$.

For any $G$-space $X$, we have the fibration
$$
X \hookrightarrow X_G \stackrel{p}{\longrightarrow} BG.
$$
The projection $p$ induces $p^*: H_G^*(pt;R) \to H_G^*(X;R)$, making $H_G^*(X;R)$ an $H_G^*(pt;R)$-module.  Natural maps in equivariant cohomology preserve this module structure.

\medskip
\noindent{\bf Borel localization}.  Let $M$ be a manifold equipped with the action of a compact torus $T$.  The first topological localization result concerns the inclusion
$$
i:M^T\hookrightarrow M
$$
of the fixed points of the action, and the map $i^*$ that this induces in equivariant cohomology.

\begin{thml1}[Borel \cite{Bo:transf}]\label{th:loc1}
Let $T$ be a compact torus, and $M$ a compact $T$-manifold.  Consider the inclusion $i: M^T\hookrightarrow M$ of the $T$-fixed points.  The kernel and the cokernel of the induced map
$$
i^*:H_T^*(M;R) \longrightarrow H_T^*(M^T;R)
$$
are {\bf torsion} $H_T^*(pt;R)$-modules.
\end{thml1}

\noindent This theorem does yield some information about the $H_T^*(pt;R)$-module structure of $H_T^*(M;R)$.  For example, if $T$ acts on $M$ freely, then $M^T=\emptyset$, so $H_T^*(M^T;R)=\{ 0\}$, and we may conclude that $H_T^*(M;R) = H^*(M/T;R)$ is entirely torsion, as an $H_T^*(pt;R)$-module.  

At the other extreme are Hamiltonian $T$-manifolds $M$:  when $R$ is a field of characteristic $0$, one can prove that $H_T^*(M;R)$ is a {\bf free} $H_T^*(pt;R)$-module.  So in this case, $i^*$ is automatically an injection.  One may also show that there is an isomorphism
$$
H^*(M;R) \cong H_T^*(M;R)\otimes_{H_T^*(pt;R)} R.
$$
This allows us to translate information about equivariant cohomology, and its ring structure, to information about ordinary cohomology.
These statements hold, with some restrictions, over $\Z$ as well.

\medskip
\noindent{\bf Atiyah-Bott Berline-Vergne localization}.  The second localization theorem relates the integral, or push-forward to a point, of an equivariant cohomology class on the whole manifold $M$ to the sum of integrals of the class over each component of the fixed point set, corrected with an equivariant characteristic class of the normal bundle to the fixed point component.

As we are using the Borel model for equivariant cohomology, we will state this result in terms of push-forwards.  Given an equivariant map $f:N\to M$ of compact oriented $G$-manifolds, there is a push-forward map in equivariant cohomology,
$$
f_*: H_G^*(N;\Q) \to H_G^*(M;\Q).
$$
To construct this map in ordinary cohomology, we use Poincar\'e duality and the natural map in homology.  For equivariant cohomology, as we are using the infinite-dimensional model $(M\times EG)/G$, we must take (large enough) finite-dimensional approximations to $EG$ to define the push-forward in any particular degree.  For a more complete discussion of the push-forward in equivariant cohomology, we refer the reader to  \cite[\S2]{ab} and \cite[Appendix~C]{ggk}.

Returning to our Hamiltonian $T$-space $M$, let $\pi : M\to pt$ be the constant map from $M$ to a point, and $\pi^F: F\to pt$ the same for a connected component of the fixed set $M^T$.  In this case of push-forward to a point, we may think of $\pi_*$ as an equivariant integral, using the Cartan model of differential forms (see, for example, \cite[Appendix~C \S 6]{ggk}). 

Just as we may define characteristic classes of vector bundles \cite{MS}, so may we associate to an {\bf equivariant} vector bundle {\bf equivariant characteristic classes}, which live in the equivariant cohomology of the base space.  The bundle in question is the normal bundle to a component  $F$ of the fixed point set $M^T$, denoted $$\nu(F\subseteq M).$$  This is a subbundle of the tangent bundle to $M$, restricted to $F$; that is,
$$\nu(F\subseteq M)\subseteq TM|_F.$$
Because $F$ is fixed, the torus $T$ acts on $TM|_F$, fixing the tangent bundle $TF$ to $F$ and acting nontrivially on the orthogonal complement  $\nu(F\subseteq M)$.  In this case, the relevant characteristic class is the {\bf equivariant Euler class} of this normal bundle, denoted $e_T(\nu(F\subseteq M))\in H_T^*(F)$.  

For further details on equivariant characteristic classes, we refer the reader to \cite[Appendix C \S6]{ggk} and references contained therein.  For our purposes here, it is sufficient to note that the classes $e_T(\nu(F\subseteq M))$ are often very easy to compute.  We may now state the push-forward version of localization.

\begin{thml2}[Atiyah-Bott \cite{ab} Berline-Vergne \cite{bv}]\label{thm:abbv}
Suppose a compact torus $T$ acts on a compact manifold $M$.  Then for any class $u\in H_T^*(M;\C)$,
\begin{equation}\label{eq:abbv}
\pi_*(u) = \sum_{{F\subseteq M^T}} \pi^F_*\left( \frac{ u|_F}{e_T(\nu(F\subseteq M))}\right) ,
\end{equation}
where the sum on the right-hand side is taken over connected components $F$ of the fixed point set $M^T$, and $u|_F$ is the restriction of $u$ to $F$.
\end{thml2}

We note that such a formula could never hold in ordinary cohomology, as ordinary Euler classes are never invertible.
This formula \eqref{eq:abbv} simplifies greatly in the case when $M^T$ consists of isolated points.  In this case, the normal bundle $\nu(F\subseteq M)$ is just the tangent space at the fixed point:
$$
\nu(F\subseteq M) = TM|_F = T_FM.
$$
One can show that in this special case, the equivariant Euler class is the product of the isotropy weights for the $T$-action on this tangent space $T_FM$. 

Finally, this localization theorem can be interpreted as a generalization of the celebrated Duistermaat-Heckman theorem \cite{dh} on the push-forward of the equivariant symplectic form.  For a discussion of this, see \cite[\S 7]{ab}.

\medskip
\noindent{\bf Localization and cobordism}.  The third localization result  that we present is the ultimate topological version of this principle.  The Atiyah-Bott Berline-Vergne localization implies that integrals of equivariant cohomology classes over $M$, for a Hamiltonian $T$-manifold $M$, can be written in terms of fixed point data.  Guillemin, Ginzburg and Karshon proved an even stronger result: that the {\bf (equivariant) cobordism class} of $M$ is in fact determined by the fixed point data.

\begin{thml3}[Guillemin-Ginzburg-Karshon \cite{ggk}]
Let $T$ be a compact torus, and suppose $M$ is a compact stably complex Hamiltonian $T$-manifold.  Then there is an equivariant cobordism,
\begin{equation}\label{eq:cobordism}
M\sim \coprod_{{F\subseteq M^T}} \nu(F\subseteq M^T)
\end{equation}
where the union on the right-hand side is taken over connected components $F$ of the fixed point set $M^T$.  In particular, if $M^T$ consists of isolated points, then $M$ is equivariantly cobordant to the disjoint union of the tangent spaces $T_pM$ for $p\in M^T$.
\end{thml3}

One must be careful to interpret \eqref{eq:cobordism} correctly: we must either compactify $\nu(F\subseteq M^T)$ via a symplectic cut, or use the notion of  {\bf proper cobordism}.  This theorem allows us to translate symplectic geometry (on the right-hand side) into representation theory (on the left-hand side).  For additional details, see \cite[Part 1]{ggk} and references contained therein.

\section{Using localization to compute equivariant cohomology}\label{sec:compute}

We return now to the first version of localization, Borel's statement that the kernel and cokernel of the map
$$
i^*:H_T^*(M;R) \hookrightarrow H_T^*(M^T;R)
$$
are torsion $H_T^*(pt;R)$-modules.  This map fits in to a long exact sequence in equivariant cohomology.  For $T= (S^1)^n$, let
$$
M_i = \{ x\in M\ |\ \dim(T\cdot x)\leq i\}
$$  
denote the {\bf $i$-skeleton} of the torus action.  In particular, then, we have $M=M_n$ and $M^T=M_0$.  These skeleta filter $M$, and there is a sequence in equivariant cohomology
\begin{equation}\label{eq:exact}
0\to H_T^*(M) \to H_T^*(M_0) \to H_T^{*+1}(M_1,M_0) \to \cdots \to H_T^{*+n}(M_n,M_{n-1}).
\end{equation}
The first (non-zero) map is precisely the map $i^*:H_T^*(M) \to H_T^*(M_0) $, and the remaining maps are defined as boundary maps in the long exact sequence of a pair. The exactness of \eqref{eq:exact} is closely related to the Serre spectral sequence  associated to the fibration
$$
\begin{array}{c}
\xymatrix{
M_{\phantom{)}}  \ar@<-0.1ex>@{^{(}->}[r] & (M\times ET)/T \ar[d] \\
& BT
}
\end{array}
$$
and whether it collapses at the $E_2$-term.  This has been studied in a general context by Franz and Puppe \cite{fp}.

An advantage of working in the symplectic category, as  Ginzburg \cite{ginz} and Kirwan \cite{K:quotients} have shown, is that for a compact Hamiltonian $T$-space, this Serre spectral sequence (over $\Q$) does indeed collapse at the $E_2$ term. 
 In particular, this implies the aforementioned fact that $H_T^*(M;\Q)$ is a {\bf free} $H_T^*(pt;\Q)$-module, and therefore that the map $i^*$ is {\bf injective}.   
Our strategy for computing $H_T^*(M;\Q)$ now will be to understand its image in $H_T^*(M^T;\Q)$.  In \cite{TW:hamTsp}, Tolman and Weitsman have shown that for a compact Hamiltonian $T$-space, the equivariant cohomology  of $M$ is isomorphic to that of its one-skeleton.  More precisely, we have

\begin{theorem}[Tolman-Weitsman \cite{TW:hamTsp}] \label{thm:oneskel}
The natural inclusions $i:M^T\to M$ and $j: M^T\to M_1$ induce the following maps in equivariant cohomology
$$
\begin{array}{c}
\xymatrix{
H_T^*(M_1;\Q) \ar[d]_{j^*} & H_T^*(M;\Q) \ar[l] \ar[dl]^{i^*} \\
H_T^*(M^T;\Q) &
}
\end{array}. 
$$
When $M$ is a compact Hamiltonian $T$-manifold, the maps $i^*$ and $j^*$ have the same image (over $\Q$).
\end{theorem}

\noindent This theorem was previously known in a more general context through work of Chang and Skjelbred \cite{CS:injectivity}.  It is also equivalent to the exactness of \eqref{eq:exact} at the term $H_T^*(M_0)$.

In order to describe the image of $i^*$ explicitly, we now make a pair of simplifying assumptions: one about the fixed point set and the other about the one-skeleton.

\begin{assum}\label{assum:fixed-points}
The fixed point set $M^T$ consists of isolated points.
\end{assum}

The point of this assumption is to simplify $H_T^*(M^T;\Q)$.  When the fixed point set consists of isolated points, this ring is a direct product of copies of
$$
H_T^*(pt;\Q)\cong \Q[x_1,\dots,x_n],
$$
one for each fixed point. Thus, every class can be represented as a tuple of polynomials, and the ring structure is the component-wise product of polynomials.

There is a large class of examples that satisfy this assumption, including symplectic toric manifolds, flag varieties, and the Hilbert scheme of $n$ points in the (complex projective) plane.  Nevertheless, this places serious restrictions on the topology of $M$.  As discussed above, a generic component $\phi^\xi$ of the moment map is a Morse function on $M$ with critical set $M^T$.  Moreover, at any fixed point, one can show that the subspace on which the Hessian of $\phi^\xi$ is negative definite  is a {\bf symplectic} subspace of the tangent space at that fixed point.  Thus, the index of each critical point is even, so $M$ is homotopy equivalent to a cell complex with only even-dimensional cells. In particular, $M$ must have trivial fundamental group.  On the other hand, Gompf has shown that any finitely presented group $G$ may be realized as the fundamental group of some compact symplectic $4$-manifold \cite{Gompf}, demonstrating just how strong this assumption is.  Karshon has studied Hamiltonian $S^1$-$4$-manifolds that satisfy Assumption~\ref{assum:fixed-points}.

\begin{theorem}[Karshon \cite{karshon}]
Let $M$ be a compact symplectic $4$-manifold equipped with a Hamiltonian $S^1$-action, satisfying Assumption~\ref{assum:fixed-points}. Then in fact $M$ is a toric surface; i.e.\  the Hamiltonian circle action extends to an effective Hamiltonian $T^2$ action.
\end{theorem}

We say that an action is {\bf semifree} if it is free outside of the fixed points.  Tolman and Weitsman have shown that Assumption~\ref{assum:fixed-points} together with the semifree assumption puts a very strong constraint on the equivariant topology of $M$.

\begin{theorem}[Tolman-Weitsman \cite{TW:semifree}]
Let $M^{2n}$ be a compact symplectic manifold equipped with a semifree Hamiltonian $S^1$ action.  If the fixed points $M^{S^1}$ are isolated, then $M$ has the cohomology and Chern classes of $(\C P^1)^k$.
\end{theorem}

\noindent In particular, a compact semifree Hamiltonian $S^1$-manifold with isolated fixed points must have exactly $2^k$ fixed points, for some $k$.    In the same paper, Tolman and Weitsman also showed that if a semifree {\em symplectic} circle action satisfies Assumption~\ref{assum:fixed-points}, then it is automatically Hamiltonian.  This is in stark contrast to McDuff's Theorem~\ref{thm:mcduff}(a).

At each $T$-fixed point $p\in M^T$, the torus acts on the tangent space to $p$.  We may always choose a $T$-invariant almost complex structure on $M$, making the tangent space $T_pM$ a complex $T$ representation.  Thus, it breaks up into one-dimensional representations,
$$
T_pM \cong \C_{\alpha_1}\oplus\cdots\oplus \C_{\alpha_d},
$$
where $\alpha_i\in \algt^*$ are the {\bf isotropy weights} of this action.  
We may phrase Assumption~\ref{assum:fixed-points} in terms of the $\alpha_i$: we simply require that all the isotropy weights at each fixed point be non-zero.  Our second assumption, which may be made independently from the first, will be stated solely in terms of these isotropy weights.

\begin{assum}\label{assum:weights}
For every point $p\in M^T$, the non-zero isotropy weights $\alpha_{i_1},\dots,\alpha_{i_k}$ are pairwise linearly independent.
\end{assum}

When this assumption is made together with Assumption~\ref{assum:fixed-points}, it is equivalent to insisting that the one-skeleton $M_1$ be two-dimensional.  We demonstrate the strength of Assumption~\ref{assum:weights} with the following theorem.

\begin{theorem}[Guillemin-Holm \cite{GuH}]
Let $M$ be a compact connected Hamiltonian $T$-space satisfying Assumption~\ref{assum:weights}.   Then all of the connected components of $M^T$ are diffeomorphic to one another.
\end{theorem}

\noindent We shall see some additional consequences of this assumption shortly.

Now combining Assumptions~\ref{assum:fixed-points} and \ref{assum:weights}, we gain a solid understanding of the equivariant topology of $M$. That the one-skeleton $M_1$ is two-dimensional is enough to imply that it is a family of embedded $\C P^1$'s, each rotated about its axis by the the torus, and intersecting one another at the fixed points.  The image of the one-skeleton under the moment map is a graph $\Phi(M_1) = \Gamma$ whose vertices correspond to the fixed points $M^T$ and whose edges correspond to the embedded $\C P^1$'s.  Each edge $e$ in $\Gamma$ is labeled by the weight\footnote{This is well-defined up to a sign, which is sufficient for our purposes.} $\alpha_e\in\algt^*$ by which $T$ acts on $e$.  Indeed, the moment map maps the corresponding $\C P^1$ to a line segment parallel to the weight $\alpha_e$. The embedding of the graph $\Gamma$ encodes, in this way, the {\bf isotropy data}, denoted $\alpha$.

The ring $H_T^*(pt;\R)$ is isomorphic to the symmetric algebra on the dual of the Lie algebra, $H_T^*(pt;\R)\cong S(\algt^*)$.  We use this to identify a torus weight $\alpha_e\in\algt^*$ with an element of $H_T^2(pt;\R)$,  the linear polynomials in the generators $x_1,\dots,x_n$.  In our examples the weights $\alpha_e$ are actually rational expressions in the $x_i$, and so we may work over $\Q$.   We now have the technical tools to state the main theorem of this section.

\begin{theorem}[Goresky-Kottwitz-MacPherson \cite{GKM}]\label{thm:gkm}
Suppose $M$ is a compact Hamiltonian $T$-space satisfying Assumptions~\ref{assum:fixed-points} and \ref{assum:weights}.  Then the image 
$$
i^*(H_T^*(M))\subseteq H_T^*(M^T)\cong \bigoplus_{p\in M^T} H_T^*(pt;\Q)
$$ consists of
\begin{equation}\label{eq:gkm}
 \left\{ (f_p)\in \bigoplus_{p\in M^T} H_T^*(pt;\Q)\ \Bigg|\ \alpha_e \big| (f_p-f_q) \mbox{ for each edge } e=(p,q) \mbox{ in } \Gamma\right\}.
\end{equation}
These divisibility conditions are often referred to as the {\bf GKM conditions}.  This ring can be defined purely from the combinatorial data $(\Gamma,\alpha)$, so is denoted $H^*(\Gamma,\alpha)$.
\end{theorem}

\begin{figure}[h]
  \begin{center}
\psfrag{x}{$x$}
  \psfrag{y}{$y$}
  \psfrag{2y}{$2y$}
  \psfrag{x+y}{$x+y$}
  \psfrag{x+2y}{$x+2y$}
 \psfrag{(a)}{(a)}
 \psfrag{(b)}{(b)}
      \epsfig{file=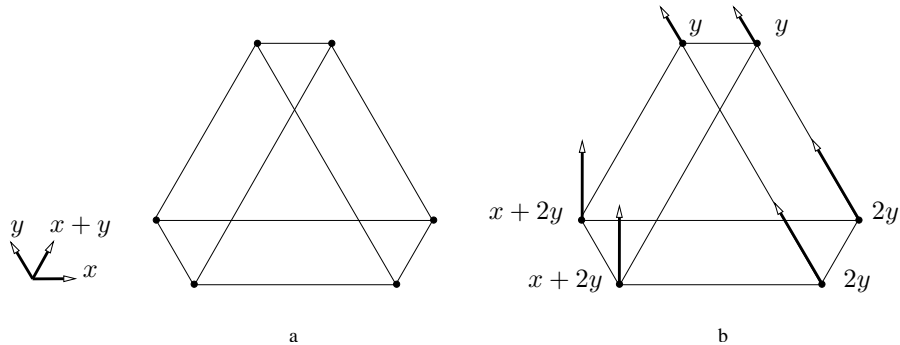,height=4.5cm}
    \caption{This figure shows (a) the moment graph for the flag variety $M = SU(3)/T = \mathcal{F}\ell(\C^3)$, with isotropy weights indicated; and  (b) a class in $H_T^*(M^T)$ that satisfies the GKM conditions. Additional details about pictorial representations of equivariant cohomology classes, such as in (b), may be found in \cite{HHH}. }\label{fig:gkm}
  \end{center}
 \end{figure}

This theorem provides a very simple prescription for computing the equivariant cohomology of a $T$-manifold satisfying Assumptions~\ref{assum:fixed-points} and \ref{assum:weights}.  We illustrate this in Figure~\ref{fig:gkm}.
Tymoczko has written an introductory account of this theorem, including many examples \cite{Ty}; for further details, the inquisitive reader could start there.
Theorem~\ref{thm:gkm} has been generalized to a wide variety of contexts. 
We conclude this section with a brief survey of some of those results.  We apologize for any (or the many!) references that may have been inadvertently excluded.

\begin{enumerate}
\item[A.] Whereas Goresky, Kottwitz and MacPherson take an algebraic approach to Theorem~\ref{thm:gkm}, Tolman and Weitsman use the Morse theory of the moment map to give an alternative proof in the symplectic category \cite{TW:hamTsp}.

\item[B.] Brion has proved an analogue for the equivariant Chow groups for nonsingular projective varieties with (algebraic) torus actions \cite{brion}.

\item[C.] Knutson and Rosu have established the result in rational equivariant K-theory  \cite[Appendix]{kr}.  Harada and Landweber have used the work mentioned in Item~E to give a careful proof for the {\em integral} K-theory of Hamiltonian $T$-spaces satisfying Assumptions $1$ and $2$ \cite{HL}.

\item[D.] Braden and MacPherson studied singular projective varieties with torus actions, and their equivariant intersection cohomology \cite{braden-macpherson}.  

\item[E.] Harada, Henriques, and the author have considered more general equivariant cohomology theories, including equivariant cobordism, in \cite{HHH}.  The setting is a bit more general; they have shown that the theory applies to a wide variety of spaces, including (infinite dimensional) coadjoint orbits of Kac-Moody groups.  They also include a discussion of coefficients.  For the particular case of equivariant $K$-theory of coadjoint orbits of Kac-Moody groups, Kostant and Kumar established the result nearly twenty years ago \cite{KK}.

\item[F.] Theorem~\ref{thm:gkm} has been extended to {\bf real loci} by Biss, Guillemin and the author \cite{BGH} and independently by Schmid \cite{schmid}.  Real loci are the symplectic analogue of the real points of a complex variety.  These results were advanced to a wider class of real loci in \cite{HH}.

\item[G.] Goldin and the author gave a combinatorial description of $H_T^*(M;\Q)$ in the case when Assumption~2 is relaxed to $\dim(M_1)\leq 4$ \cite{GoH}.  In this case, the Atiyah-Bott Berline-Vergne Localization Theorem~\ref{thm:abbv} plays a crucial role.

\item[H.] Braden, Chen and Sottile have computed the equivariant cohomology of {\bf quot schemes} \cite{BCS}. These spaces satisfy Assumption~1 but not Assumption~2. The authors explicitly compute the equivariant cohomology of each of the pieces of the one-skeleton.  These pieces all turn out to be toric varieties for various extensions of the torus.

\item[I.] Guillemin and the author worked on Hamiltonian $T$-manifolds that satisfy Assumption 2 only.  This is sufficient to ensure that $\Phi(M_1)$ is still a graph, where the vertices  now each correspond to a  connected component of the fixed set $M^T$. It turns out that all components of $M^T$ are diffeomorphic to the same manifold, denoted $F$.  The main theorem of \cite{GuH} is that under Assumption 2,
$$
H_T^*(M;\Q) \cong H^*(F)\otimes H^*(\Gamma,\alpha).
$$

\item[J.] McMullen has exploited the relationship between symplectic toric manifolds and convex polytopes to give a new combinatorial proof of the hard Lefshetz  theorem for symplectic toric manifolds \cite{McMullen}.  While this theorem is true more generally, it is a deep and difficult result.  McMullen's work provides nice insight into the workings of this theorem in the specific case of  symplectic toric manifolds.

\end{enumerate}

\section{Combinatorial localization and polytope decompositions}\label{sec:polytopes}

We conclude this article with a very combinatorial localization phenomenon.  Motivated by the fact that the vertices of a moment polytope correspond to the fixed points of the Hamiltonian torus action, we decompose an arbitrary polytope $\Delta$ in $\R^n$  (equipped with the usual inner product) into an alternating sum of cones, one for each fixed point.  For each vertex $v$ in $\Delta$, let $E_v = \{ \alpha_{v,1},\dots,\alpha_{v,k_v}\}$ be the set of edge vectors emanating from $v$.  Each $\alpha_{v,j}$ is determined up to a positive scalar.  In terms of these edge vectors, the {\bf tangent cone} at a vertex $v$ is the positive span of these edges:
$$
C_v := \left\{ v +\sum_{\alpha_{v,j}\in E_v} x_j\alpha_{v,j} \  \ \Bigg| \  x_j\geq 0 \mbox{ for all } j\right\}.
$$
This is a polyhedral cone.  The polar decomposition theorem expresses the characteristic function of the polytope as a linear combination of the characteristic functions of polarizations of these tangent cones.  To define the {\em polarizations}, we fix a vector $\xi\in \R^n$ satisfying $\alpha_{v,j}\cdot \xi \neq 0$ for all $v$ and $j$.  We call such a $\xi$ a polarizing vector, and think of it as defining ``upwards" in $\R^n$.  

First we polarize the edge vectors so that they all point ``downwards".  Let
$$
\alpha_{v,j}^\# = \left\{ \begin{array}{rl}
	\alpha_{v,j} & \mbox{ if } \alpha_{v,j}\cdot \xi < 0, \mbox{ and}\\
	-\alpha_{v,j} & \mbox{ if } \alpha_{v,j}\cdot \xi >0.
	\end{array}\right. 
$$
We keep track of the sign changes by letting
$$
E_v^+ := \left\{ \alpha_{v,j}\ | \  \alpha_{v,j}\cdot \xi >0\right\} \mbox{ and }
E_v^- := \left\{ \alpha_{v,j}\ | \  \alpha_{v,j}\cdot \xi <0\right\}.
$$
We now define the {\bf polarized tangent cone} at a vertex $v$ to be
$$
C_v^\# := \left\{  \left. v +\sum_{\alpha_{v,j}\in E_v} x_j\alpha_{v,j}^\#\  \right| \begin{array}{rl}
	x_j\geq 0 & \mbox{ if } \alpha_{v,j}\in E_v^-, \mbox{ and}\\
	x_j>0 & \mbox{ if } \alpha_{v,j}\in E_v^+
\end{array} \right\}.
$$
Finally, recall that the characteristic function of a set $A\subseteq\R^n$ is
$$
{\mathbf{1}}_A(x) = \left\{ \begin{array}{ll}
1 & x\in A,\\
0 & x\not\in A.
\end{array}\right.
$$

\begin{theorem}[Lawrence \cite{Lawrence}, Varchenko \cite{Varchenko}]
For any convex polytope $\Delta\subseteq \R^n$,
$$
{\mathbf{1}}_\Delta(x) = \sum_{v\in\Delta} (-1)^{\left|E_v^+\right|} {\mathbf{1}}_{C_v^\#}(x).
$$
\end{theorem}

\begin{figure}[h]
  \begin{center}
  \psfrag{=}{$=$}
  \psfrag{-}{$-$}
  \psfrag{+}{$+$}
  \psfrag{x}{$\xi$}
      \epsfig{file=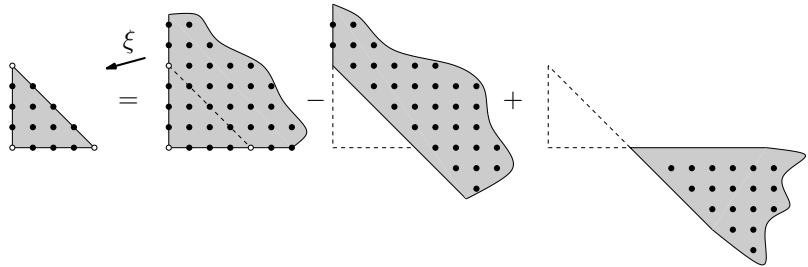,height=3.5cm}
    \caption{The polar decomposition theorem for the triangle, for the indicated value of $\xi$. }
    \label{fig:polar}
  \end{center}
 \end{figure}

We illustrate the result in a simple example in Figure~\ref{fig:polar}.
Karshon, Sternberg and Weitsman give a short direct proof of this theorem in \cite{ksw2}. They use this to give a combinatorial proof of the Euler-Maclaurin formula that relates the sum of values of a (nice) function $f$ at the lattice points inside a polytope on the one hand, to the integral of $f$ over the polytope on the other.  While the technical details of this work are beyond the scope of this article, the curious reader can find particulars in \cite{ksw1,ksw2}.

\end{document}